\documentclass[12pt]{article}
\usepackage{amssymb}
\usepackage{amsmath, amsthm}
\usepackage{multirow}
\usepackage[normalem]{ulem} 

\usepackage{enumitem} 
\setlist[enumerate]{topsep=2pt,itemsep=-0ex,partopsep=1ex,parsep=1ex} 

\usepackage{hyperref}
\hypersetup{
  colorlinks   = true, 
  urlcolor     = blue, 
  linkcolor    = blue, 
  citecolor   = red 
}

\newcommand{\G}{\Gamma}
\newcommand{\cstar}{\bb{C}^\times}

\newcommand{\diag}[1]{\mbox{diag}\left\{#1\right\}}

\newcommand{\eqr}[1]{\eqref{eq:#1}}

\newcommand{\Hom}{\mbox{Hom}}
\newcommand{\e}[1]{\textbf{e}\left(\textstyle#1\right)}
\newcommand{\mc}[1]{\mathcal{#1}}

\newcommand{\cvec}[3]{\begin{pmatrix}#1\\#2\\#3\end{pmatrix}}
\newcommand{\bb}[1]{\mathbb{#1}}
\newcommand{\gln}[1]{\mbox{GL}_{#1}(\bb{C})}

\newcommand{\C}{\bb{C}}

\newcommand{\up}{\upsilon}
\newcommand{\tr}{\mbox{tr}}   
\newcommand{\wt}{\mbox{wt}\,}

\newtheorem{thm}{Theorem}[section]
\newtheorem{lem}[thm]{Lemma}

\newtheorem{cor}[thm]{Corollary}
\newtheorem{prop}[thm]{Proposition}

\newtheorem{rmk}[thm]{Remark}

\numberwithin{equation}{section}

\begin{document}

\title{Intertwining operators and vector-valued modular forms for minimal models}
\author{\begin{tabular}{c}
Matthew Krauel\footnote{e-mail: krauel@csus.edu. Portions of this research and preliminary work for this paper had been supported by the European Research Council (ERC) Grant agreement n. 335220 - AQSER, as well as the Japan Society of the Promotion of Science (JSPS), No. P13013} \\
{\small Department of Mathematics and Statistics, California State University, Sacramento, USA} \\ 
{\small and} \\
Christopher Marks\footnote{e-mail: cmarks@csuchico.edu} \\
 {\small Department of Mathematics and Statistics, California State University, Chico, USA}
\end{tabular}}
\date{}

\maketitle
\begin{abstract}
Using the language of vertex operator algebras (VOAs) and vector-valued modular forms we study the modular group representations and spaces of $1$-point functions associated to intertwining operators for Virasoro minimal model VOAs. We examine all representations of dimension less than four associated to irreducible modules for minimal models, and determine when the kernel of these representations is a congruence or noncongruence subgroup of the modular group. Arithmetic criteria are given on the indexing of the irreducible modules for minimal models that imply the associated modular group representation has a noncongruence kernel, independent of the dimension of the representation. The algebraic structure of the spaces of 1-point functions for intertwining operators is also studied, via a comparison with the associated spaces of holomorphic vector-valued modular forms.
\end{abstract}

\section{Introduction\label{Sect:Intro}}

The motivation behind this paper lies in a desire by the authors to better understand modularity in the setting of intertwining operators for rational vertex operator algebras (VOAs). As we elaborate on momentarily, we seek to understand both the arithmetic of the Fourier coefficients of 1-point functions arising in this setting, as well as the algebraic structure of the graded spaces of vector-valued modular forms into which these 1-point functions organize themselves.
The first author learned from Masahiko Miyamoto that one-dimensional examples of vector-valued 1-point functions had already been produced by him \cite{Miy-Int} for certain Virasoro minimal models. Seeking to build on these examples, a subsequent examination of the literature revealed that  nearly nothing had been published regarding modular group representations arising from intertwining operators for rational VOAs, not to mention the arithmetic and algebraic considerations of the spaces of vector-valued modular forms arising in this setting (a notable exception is the work of Etinghof and Kirillov \cite{EK-Jack, Kirillov-Inner}).

Given this state of affairs, the authors decided it would be useful to initiate a general study of modularity for intertwining operators associated to rational VOAs. We quickly realized, however, that it was both natural and practical to focus our attention on the Virasoro minimal models, as this provides a somewhat elementary setting in which to work as well as serving as a direct generalization of the examples found in \cite{Miy-Int}. As we demonstrate below, for the questions we are interested in we are able to find more or less complete answers when studying modular group representations of dimension less than four arising from intertwining operators for Virasoro modules. These investigations form the heart of the paper.

On the arithmetic side, the second author has been interested for some time in noncongruence modular forms, and was made aware (via conversations with Terry Gannon) that intertwining operators may be viewed as machines for producing copious examples of noncongruence vector-valued modular forms. (We recall here that a vector-valued modular form is called noncongruence if there is no positive integer $N$ such that the kernel of the representation according to which the vector transforms contains the principal congruence subgroup
\begin{equation}\label{eq:gammaN}
\G(N)=\left\{\begin{pmatrix}a&b\\c&d\end{pmatrix}\in \operatorname{SL}_2(\bb{Z})\vert a\equiv d\equiv1, b\equiv c\equiv0\pmod{N}\right\}
\end{equation}
of level $N$.) In contrast to vertex operators for rational VOAs, where it has long been expected (see \cite{DLinNg} and references therein for results and a review of relevant literature) that kernels of associated modular group representations must be congruence, it is believed in the more general setting of intertwining operators that these representations are generically noncongruence. There is ample evidence that the presence (or lack) of a congruence subgroup in the kernel of such a representation has a profound impact on the arithmetic of the Fourier coefficients of the 1-point functions associated to the representation. Indeed, a longstanding conjecture originating in an article \cite{ASD} of Atkin and Swinnerton-Dyer (and subsequently extended to the setting of vector-valued modular forms by Geoffrey Mason \cite{Mason}) predicts that only congruence modular forms may have Fourier coefficients contained in the ring of integers of a number field. In particular, if a given vector-valued modular form arising in the setting of intertwining operators for a rational VOA were to have a combinatorial interpretation -- by which we mean that its Fourier coefficients are rational integers (perhaps even natural numbers) that count some data associated with the VOA or its irreducible modules -- then according to the above-mentioned Bounded Denominator Conjecture the associated modular group representation is congruence (i.e., has a congruence subgroup as kernel). Thus the above conjecture motivates an identification of the congruence representations arising from intertwining operators in rational VOA theory.

In Section \ref{sec:vvmf} below we investigate the nature of the modular group representations arising from intertwining operators for Virasoro minimal models, both in the low-dimension setting and in arbitrary dimension. In particular, in Theorems \ref{thm:dim1}, \ref{thm:dim2}, and \ref{thm:dim3} below we give a complete classification of such representations when the dimension is one, two, or three respectively. Already one sees in this setting that all three possibilities occur, i.e.,  the kernel of the representation could be a congruence subgroup, a finite index noncongruence subgroup, or an infinite index subgroup (thus noncongruence). We also prove some arithmetic results that apply in arbitrary dimension. By utilizing work of Nobs and Wolfart \cite{Nobs, Nobs-Wolfart} concerning representations of $\operatorname{SL}_2(\bb{Z}/N\bb{Z})$ we are able to give sufficiency conditions -- independent of the dimension of the representation -- that imply it has a noncongruence kernel. Explicitly, we study $V(p,q)$ when $p$ and $q$ are prime powers and show in Theorem \ref{thm:pq} below that for such a minimal model most irreducible $V(p,q)$-modules yield a noncongruence representation of the modular group. A particularly clean result follows from the proof of Theorem \ref{thm:pq}, which we state here to give a flavor of what is found in Section \ref{sec:vvmf}. As discussed there, we employ the standard Kac indexing for the irreducible $V(p,q)$-modules, writing them as $L_{m,n}$ for pairs $(m,n)$ of odd integers in the ranges $1\leq m\leq p-1$, $1\leq n\leq q-1$.
\begin{cor}\label{cor:pq-2}
Suppose $L_{m,n}$ is an irreducible $V(p,q)$-module with $p,q>3$ distinct primes and 
$$(m,n) \not \in \{(1,1),(1,q-2),(p-2,1),(p-2,q-2)\}.$$ 
Then the modular group representation $\rho_{m,n}$ arising from the action of $L_{m,n}$ is noncongruence.\qed
\end{cor}

We note that this result utilizes our conventions \eqr{conventions} below regarding the parity of $p$, $m$, and $n$. We also point out that when $(m,n)=(1,1)$ then $L_{1,1}=V(p,q)$ is the vacuum module so that $\rho_{m,n}=\rho_{1,1}$ is known to be congruence \cite{DLinNg, Rocha}, and similarly when $(m,n)=(p-2,q-2)$ it follows that $\rho_{p-2,q-2}$ is a two-dimensional congruence representation covered by case $(i)$ of Theorem \ref{thm:dim2} below.

Beyond the above arithmetic considerations, we are also interested in the algebraic structure of spaces of 1-point functions associated to intertwining operators. As we discuss in Section \ref{sec:vvmf} below, work of Miyamoto \cite{Miy-Int} (also Huang \cite{Huang-Int} for the $n$-point case, and Yamauchi \cite{Yamauchi-Int} for the orbifold case) implies that for a rational VOA $V$ and an irreducible $V$-module $U$ the 1-point functions arising from the intertwining operators associated to $U$ may be organized into a graded space $\mc{V}(\rho)$ of weakly holomorphic vector-valued modular forms, where $\rho$ denotes the modular group representation carrying the action of the intertwining operators for $U$. It follows from Lemma \ref{lem:surjective1} below that $\mc{V}(\rho)$ contains a substantial portion of the related space $\mc{H}(\rho)$ of holomorphic vector-valued modular forms for $\rho$, and we are interested in knowing when, in fact, these two spaces coincide. This generalizes a question asked by Dong and Mason in \cite{DM-Higher}, where it was determined which modular forms for the full modular group are realized as 1-point functions for the Moonshine module. Again we are able to give a complete answer to this question in the low-dimension setting of Virasoro minimal models, and in Theorems \ref{thm:dim1}, \ref{thm:dim2}, and \ref{thm:dim3} below we give bounds on the ratio $\frac{q}{p}$ that imply equality between the spaces $\mc{V}(\rho)$ and $\mc{H}(\rho)$ when $\dim\rho<4$ (it also follows that when these bounds are exceeded there will be proper containment $\mc{V}(\rho)\subset\mc{H}(\rho)$).

The authors would like to thank Terry Gannon (who first made us aware of the existence of noncongruence representations in VOA theory) and Masahiko Miyamoto (whose work we use in a crucial manner here) for many fruitful conversations.


\section{Preliminaries\label{Sect:Preliminaries}}
Let $V$ be a simple rational VOA with central charge $\mathbf{c}\in\bb{Q}$. Then $V$ has finitely many inequivalent irreducible modules \cite{DLM-Twisted}, which we denote here as $V=W^1,\ldots,W^d$. Let $h_i\in\bb{Q}$ denote the conformal weight of $W^i$ and $Y^i(u ,z):=\sum_{n\in \mathbb{Z}} u(n)z^{-n-1}$ be the vertex operator associated with an element $u$ of $W^i$, though we often simply write $Y$ when the dependence on $W^i$ is clear or unnecessary. Since $V$ is rational, each $W^i$ has an $L(0)$-induced grading $W^i=\bigoplus_{n\geq 0}W^i_{h_i +n}$, where the endomorphisms $L(m)$ are derived from the Virasoro element $\omega \in V$ by $Y (\omega ,z) :=\sum_{m\in \bb{Z}} \omega(m)z^{-m-1}=:\sum_{m\in \bb{Z}}L(m)z^{-m-2}$. When a homogeneous element $u$ is contained in $W^i_{\alpha}$, we say $u$ is of weight $\alpha$. If the precise weight is unknown with respect to the $L(0)$-grading, we say $u$ is of weight $\wt u$. Zhu showed in \cite{Zhu} that setting $Y[v,z]=Y^i [v ,z]:=Y^i (e^{zL(0)}v,e^z-1)=:\sum_{m\in \mathbb{Z}} v[m]z^{-m-1}$ for homogeneous $v\in V$ produces another $V$-module structure on $W^i$, and by defining $L[m]=\widetilde{\omega} [m+1]$ (where $\widetilde{\omega} :=\omega -\mathbf{c}/24$) via these vertex operators, $L[0]$ induces a different grading which we write as $W^i=\bigoplus_{n\geq 0}W^i_{[h_i +n]}$. We denote the weight of a homogeneous element $u \in W^i$ with respect to the $L[0]$-grading by $\wt[u]$.\\
\indent There is a fusion product defined for these modules in the associated modular tensor category, where for each $i$ and $j$ we have
\begin{equation}\label{eq:fusion}
W^i\times W^j=\sum_{k=1}^dN_{i,j}^kW^k.
\end{equation}
The fusion coefficient $N_{i,j}^k$ is the dimension of the $\bb{C}$-vector space $\mc{I}_{i,j}^k$ of \emph{intertwining operators of type $\binom{W^k}{W^i \, W^j}$}, which are linear maps
\begin{eqnarray*}
\mathcal{Y}(\ast,z)\colon W^i&\rightarrow&\Hom(W^j,W^k)\{z\},\\
w_i&\mapsto& \mathcal{Y}(w_i,z)=\sum_{r\in\bb{C}}w_i(r)z^{-r-1}
\end{eqnarray*}
that satisfy the following conditions for all $w_i \in W^i$, $w_j \in W^j$, and $v\in V$:
\begin{enumerate}
\item $w_i (r) w_j =0$ for $r$ sufficiently large,
\item $\mathcal{Y}(L(-1)w_i ,z) =\frac{d}{dz} \mathcal{Y}(w_i ,z)$, and
\item the Jacobi identity
\begin{align*}
 z_0^{-1} \delta &\left( \frac{z_1 -z_2}{z_0}\right) Y^k (v,z_1)\mathcal{Y}(w_i,z_2)w_j -z_1^{-1} \delta \left( \frac{z_2 -z_1}{-z_0} \right) \mathcal{Y}(w_i ,z_2) Y^j (v,z_1)w_j \\
 &=z_2^{-1} \delta \left( \frac{z_1 -z_0}{z_2} \right) \mathcal{Y}\left(Y^i (v,z_0)w_i,z_2\right)w_j.
\end{align*}
\end{enumerate}

 \indent It is a consequence of this definition and the $L(0)$-grading discussed above, that for any homogeneous $w_i\in W^i$ we have
 \begin{equation}\label{eq:N-grading}
 w_i(k)\colon W^j_{h_j+n} \to W^k_{h_j+n+\wt(w_i)-k-1}.
 \end{equation}
  In particular, in the case $W^j=W^k$ and $k=\wt(w_i)-1$, we have $w_i(k)W^j_{h_j+n}\subseteq W^j_{h_j+n}$. That is, the endomorphism $o^{\mathcal{Y}}(u):=u(\wt(u)-1)$, often called the \textit{zero mode} of $u$ (with respect to $\mathcal{Y}(\ast ,z)$), preserves all weight spaces. When it is clear which intertwining operator we are discussing, we simply write $o(u)$.

 Similar to the vertex operator structure $Y[\ast ,z]$ discussed above, we can define another intertwining operator by setting $\mathcal{Y}[w_i,z]:=\mathcal{Y}(e^{zL(0)}w_i,e^z-1 )=:\sum_{r\in \mathbb{C}}w_i[r] z^{-r-1}$ for $w_i \in W^i$. Then $\mathcal{Y}[\ast ,z]$ is an intertwining operator of type $\binom{W^k}{W^i\, W^j}$ for the modules $(W^i,Y^i[\ast])$, $(W^j,Y^j[\ast])$, and $(W^k,Y^k[\ast])$. Note that a similar grading to that of (\ref{eq:N-grading}) occurs for endomorphisms $w_i[r]$ and the $L[0]$-grading.
 
 A $V$-module $W$ is called \textit{$C_2$-cofinite} if the subspace 
\[ 
C_2 (W) :=\mbox{span}_\bb{C}\{a(-2)w \mid a \in V, w\in W\}
\] 
has finite codimension in $W$. Miyamoto \cite{Miy-Int} (and also Yamauchi \cite{Yamauchi-Int}) utilized a concept of $C_{[2,0]}$-cofinite for modules. However, \cite[Proposition 5.2]{ABD} gives that every irreducible module of a $C_2$-cofinite VOA is itself $C_2$-cofinite. Along with the fact that $C_2$-cofinite implies $C_{[2,0]}$-cofinite, we only need to be concerned that $V$ is $C_2$-cofinite to apply the results of \cite{Miy-Int, Yamauchi-Int}. We therefore assume this for the remainder of the article.

 Fix a $V$-module $U$ and let $S=S_i$ denote the subset of $\{1,2,\dots,d \}$ consisting of all $j$ such that $N_{i,j}^j>0$. For each $j\in S$, fix a basis $\{\mathcal{Y}_j^{(k)}(*,z)\mid 1\leq k \leq N_{i,j}^j \}$ of $\mc{I}_{i,j}^j$. For $u\in U$ and an intertwining operator $\mathcal{Y}_j^{(k)} \in \mc{I}_{i,j}^j$,  consider the function defined by
\begin{equation}\label{eq:tracefcn}
S_j^{(k)}(u,\tau) :=\tr|_{W^j}o(u)q^{L(0)-\frac{\mathbf{c}}{24}}=q^{h_j-\frac{\mathbf{c}}{24}}\sum_{n\geq0}\tr|_{W^j_{h_j+n}} o (u)q^n,
\end{equation}
and extend this definition linearly in $u$. 
Given the $V$-module $U$ as above, we denote by $\mathcal{C}_1 (U)$ the space of $1$-point functions for $W^i$. More precisely, 
\begin{equation}
\label{eq:1ptSpace}
\mathcal{C}_1 (U) :=\mbox{span}_\bb{C}\left\{S_j^{(k_j)} (w,\tau) \mid j\in S_i, 1\leq k_j \leq N_{i,j}^j, w\in U \right\}.
\end{equation}
Additionally, for $u\in U$ set
\begin{equation}\label{eq:tracefcnC1}
\mathcal{C}_1^u (U) :=\mbox{span}_\bb{C}\left\{S_j^{(k_j)} (u,\tau) \mid j\in S_i, 1\leq k_j \leq N_{i,j}^j \right\}.
\end{equation}
By results of Miyamoto \cite{Miy-Int} (and also Huang \cite{Huang-Int} and Yamauchi \cite{Yamauchi-Int}), each function \eqref{eq:tracefcn} is holomorphic in $\bb{H}$ and the space $\mathcal{C}_1^u (U)$ for homogeneous $u\in W^i$ with respect to $L[0]$ can be described as a weakly holomorphic vector-valued modular form of weight $\wt[u]$ for the full modular group $\G :=\mbox{SL}_2(\bb{Z})$. 

Unlike in the setting of vertex operators, where $\wt[u]$ is integral, if $U\neq V$ then generically $\wt[u]\in\mathbb{Q}\setminus \mathbb{Z}$. Consequently, one must choose a multiplier system in order to obtain a representation of $\G$. (A detailed discussion of real weight modular forms is found in \cite[Chapter 3]{Rankin}, and an extensive discussion of vector-valued modular forms for arbitrary real weight (along with additional references) is found in the second author's doctoral dissertation \cite{Marks1}.) For each weight, there are 12 possible choices of multiplier system which differ only by a character (i.e., one-dimensional representation) of $\G$. Since each such character has a congruence kernel (of level dividing 12), this choice of character has no effect on the results obtained in this article. Thus, for the sake of definiteness, we fix once and for all the following convention: for each $k\in\bb{R}$ we fix $\up=\up_k$ to be the unique multiplier system in weight $k$ such that 
\begin{equation}\label{eq:upT}
\up \begin{pmatrix} 1 & 1\\0 &1 \end{pmatrix} =\e{\frac{k}{12}},
\end{equation} 
where for $r\in\bb{R}$ we set $\e{r}=e^{2\pi ir}$.
Thus $\up$ is the multiplier system according to which $\eta^{2k}$ transforms as a modular form of weight $k$, where 
\begin{equation}\label{eq:eta}
\eta (\tau) :=q^\frac{1}{24}\prod_{n\geq1}(1-q^n)
\end{equation} 
denotes Dedekind's eta-function (here $\tau$ is a variable in the complex upper half-plane $\bb{H}$ and $q=e^{2\pi i\tau}$). Given a $d$-dimensional representation $\rho:\G\rightarrow\gln{d}$ of the modular group, we then say that a holomorphic function $F=(f_1,\ldots,f_d)^t:\bb{H}\rightarrow\bb{C}^d$ is a weakly holomorphic vector-valued modular form of weight $k$ for $\rho$ if $F$ is meromorphic at the cusps of $\G$ and satisfies the functional equation
$$F|_k^\up\gamma=\rho(\gamma)F.$$
Here we define, for each $\tau\in\bb{H}$ and $\gamma=\left( \begin{smallmatrix} a&b \\ c&d \end{smallmatrix} \right) \in\G$,
\begin{equation}\label{eq:slash}
F|_k^\up\gamma(\tau)=\up(\gamma)^{-1}(c\tau+d)^{-k}F\left(\frac{a\tau+b}{c\tau+d}\right).
\end{equation}

We return to the context of the current article. Choosing a basis  
\begin{equation}\label{eq:basis}
\left\{S_{j_1}^{(k_1)}(u,\tau),\ldots, S_{j_s}^{(k_s)}(u,\tau)\right\}
\end{equation}
for the space $\mathcal{C}_1^u(U)$ associated to an element $u\in U_{[\wt [u]]}$, one obtains a matrix representation 
\begin{equation}\label{eq:rho}
\rho:\G\rightarrow\gln{s}
\end{equation} 
according to which 
\begin{equation}\label{eq:F}
F(u,\tau)=\cvec{S_{j_1}^{(k_1)}(u,\tau)}{\vdots}{S_{j_s}^{(k_s)}(u,\tau)}=\cvec{q^{h_1-\frac{\bf{c}}{24}}\sum\limits_{n\geq0}\tr|_{W^{j_1}_{h_1+n}}o(u)q^n}{\vdots}{q^{h_s-\frac{\bf{c}}{24}}\sum\limits_{n\geq0}\tr|_{W^{j_s}_{h_s+n}}o(u)q^n}
\end{equation} 
transforms as a weakly holomorphic vector-valued modular form of weight $\wt[u]$. Note that the dimension $s$ of \eqr{rho} is bounded by the sum of the relevant fusion coefficients, i.e., the inequality
\begin{equation}\label{eq:dimbound}
s\leq\sum_{j\in S_i}N_{i,j}^j
\end{equation}
holds. Equality in \eqr{dimbound} occurs only if all the $1$-point functions are linearly independent. In general, depending on the VOA, modules, and elements involved, linear dependence among $1$-point functions can occur, thus complicating the computation of $s$.  
  
Let $\mathcal{M}:=\bigoplus_{j\geq0}\mc{M}_{2k}$ denote the graded ring of holomorphic modular forms for $\G$. Thus $\mc{M}=\mathbb{C}[G_4,G_6]$ where for even $k\geq2$ we denote\footnote{We note that the functions $G_k$ defined here are precisely those written as $E_k$ in \cite{DM-Higher}.} by $G_k$ the Eisenstein series in weight $k$
 \begin{equation}\notag
  G_k(\tau):=-\frac{ B_{k}}{k!}+\frac{2}{(k-1)!}
\sum\limits_{n\geq 1}\frac{n^{k-1}q^{n}}{1-q^{n}}
\end{equation}
(here $B_k$ denotes the $k$th Bernoulli number). 
With this notation, the \emph{modular derivative} in weight $k\in\bb{R}$ is written
\begin{align*}
\partial_k&=\frac{1}{2\pi i}\frac{d}{d\tau}+kG_2 = q\frac{d}{dq}+kG_2.
\end{align*}
For each $k\in\bb{R}$, $\partial_k$ is covariant with respect to the slash action of $\G$, i.e., for any meromorphic $f:\bb{H}\rightarrow\bb{C}$ and any $\gamma\in\G$ we have
\begin{equation}\label{eq:covar}
(\partial_kf)|_{k+2}^\up\gamma=\partial_k(f|_k^\up\gamma),
\end{equation}
and in particular if $f$ is modular of weight $k$ then $\partial_kf$ is modular of weight $k+2$. Similarly, letting $\partial_k$ act componentwise it follows that if $F$ is a vector-valued modular form of weight $k$ for a representation $\rho$ then $\partial_kF$ is of weight $k+2$ for the same representation $\rho$. This prompts one to define the ring of modular differential operators
\begin{equation}\label{eq:R}
 \mathcal{R}=\{ \phi_0 + \phi_1 \partial + \cdots + \phi_n \partial^n \mid \phi_i \in \mathcal{M}, n\geq 0 \},
 \end{equation}
where addition is performed as if $\mc{R}$ were the polynomial ring $\mc{M}[\partial]$, and multiplication in $\mc{R}$ is defined via the ``Leibniz rule'' $\partial\cdot\phi=\phi\partial+\partial_k\phi$ for any $\phi\in\mc{M}_k$. In the next section, we will be concerned with the cyclic module $\mc{R}F$ of $1$-point functions, with $F$ as in \eqr{F} (cf.\ Theorem \ref{thm:1ptgen} below).
 
Returning to our discussion regarding the space of $1$-point functions for a $V$-module $U$, recall that a primary vector $u\in U$ satisfies $L(n)u=0$ for $n\geq 1$. Then the proof of Proposition 2(b) in \cite{DM-Higher} implies the following result.

 \begin{lem}\label{lem:surjective1}
Suppose $u\in U_{[h_i]}$ is a primary vector. Then for any $j\in S_i$ and $1\leq k_j \leq N_{i,j}^j$, we have $\mathcal{R} S_j^{(k_j)}(u,\tau) \subseteq \mathcal{C}_1 (U)$. That is, for any $f(\tau) \in \mathcal{R} S_j^{(k_j)}(u,\tau)$, there exists an element $w\in U$ such that $S_j^{(k_j)}(w,\tau) = f(\tau)$. \hfill \qed
\end{lem}
 

\section{Minimal models and vector-valued modular forms\label{sec:vvmf}}

In this section we initiate a detailed study of the modular group representations and spaces of $1$-point functions arising from intertwining operators related to Virasoro minimal models. We are interested both in the nature of the representations arising in this setting and also in the structure of the space of $1$-point functions associated to such a representation. We first review  relevant definitions and results, and then commence with an analysis of the representations of dimension less than four arising from intertwining operators for irreducible modules of Virasoro minimal models. We conclude by proving some results that hold in arbitrary dimension.

The Virasoro minimal models are a family $\{V(p,q)\}$ of rational VOAs, indexed by pairs $(p,q)$ of relatively prime integers $p,q\geq2$. The VOA $V(p,q)$ is uniquely determined by its central charge  
\begin{equation}\label{eq:cc}
\mathbf{c}=\mathbf{c}_{p,q}=1-\frac{6(p-q)^2}{pq}.
\end{equation}
We employ the standard Kac indexing for the irreducible $V(p,q)$-modules, which is given by pairs of integers $(m,n)$ appearing in what is termed (we use \cite[Theorem 4.2]{Wang} as a reference) an \emph{admissible triple}
\begin{equation}\label{eq:rules}
\{(m,n),(m_j,n_j),(m_k ,n_k)\}.
\end{equation}
The rules for admissible triples are as follows:
\begin{enumerate}
\item[(A1)] $0<m,m_j,m_k<p$,\ \ $0<n,n_j,n_k<q$. 
\item[(A2)] $m<m_j+m_k$, $m_j<m+m_k$, $m_k<m+m_j$, and $n<n_j+n_k$, $n_j<n+n_k$, $n_k<n+n_j$.
\item[(A3)] $m+m_j+m_k<2p$,\ \ $n+n_j+n_k<2q$. 
\item[(A4)] $m+m_j+m_k$ and $n+n_j+n_k$ are odd. 
\item[(A5)] The identifications 
\begin{equation}\label{eq:triples}
\{(m,n),(m_j,n_j),(m_k,n_k)\}\sim\{(m,n),(p-m_j,q-n_j),(p-m_k,q-n_k)\} \notag
\end{equation}
are enforced.
\end{enumerate}
Since $p$ and $q$ are relatively prime and \eqr{cc} is symmetric with respect to $p$ and $q$, we may (and shall) assume that $p\geq3$ is odd. With this assumption and the above rules, one may index the irreducible modules of $V(p,q)$ by the integer pairs 
\[\left\{(m,n) \, \biggl \vert \, 1\leq m\leq\frac{p-1}{2},\ \ 1\leq n\leq q-1\right\}.\]
This yields $\frac{(p-1)(q-1)}{2}$ irreducible $V(p,q)$-modules. The conformal weight of the module $L_{m,n}$ corresponding to the pair $(m,n)$ is
\begin{equation}\label{eq:cw}
h=h_{m,n}=\frac{(np-mq)^2-(p-q)^2}{4pq}.
\end{equation} 
In particular, setting $m=n=1$ yields the vacuum module $L_{1,1}=V(p,q)$ with conformal weight $h_{1,1}=0$. Since the irreducible $V(p,q)$-modules are uniquely determined by their conformal weight and evidently $h_{m,n}=h_{p-m,q-n}$ for any $m,n$, we have $L_{m,n}=L_{p-m,q-n}$. We make crucial use of this fact in what follows.

Let $\mathcal{I}_{(m ,n),(m_j,n_j)}^{(m_k,n_k)}$ denote the space of intertwining operators of type $\binom{L_{m_k,n_k}}{L_{m,n}\, L_{m_j,n_j}}$, and let $N_{(m_i,n_i), (m_j,n_j)}^{(m_k,n_k)}=\dim_\bb{C}\mathcal{I}_{(m ,n),(m_j,n_j)}^{(m_k,n_k)}$ be the corresponding fusion rule. Then \cite[Theorem 4.3]{Wang} $N_{(m_i,n_i), (m_j,n_j)}^{(m_k,n_k)}=0$ unless $\{(m,n),(m_j,n_j),(m_k,n_k)\}$ is an admissible triple, and in this case $N_{(m_i,n_i), (m_j,n_j)}^{(m_k,n_k)}=1$. Since $L_{m_j,n_j}=L_{p-m_j,q-n_j}$ for any $(m_j,n_j)$, there are two basic types of admissible triples that capture the action of a fixed module $L_{m,n}$ on the other irreducible $V(p,q)$-modules, namely
\[
\{(m,n),(m_j,n_j),(m_j,n_j)\} \quad \mbox{and} \quad \{(m,n),(p-m_j,q-n_j),(m_j,n_j)\}.
\]
By rule (A4) above, we require in the above cases that either $m+2m_j$ or $m+p$ be odd, and since we are assuming that $p$ is odd this implies that $m$ is odd in the former case but even in the latter. Thus we may (and shall) adopt the further convention that when the action of the $V(p,q)$-module $L_{m,n}$ is considered, the integer $m$ will be taken to be odd. Note that since $L_{m,n}=L_{p-m,q-n}$ and $p-m$ is even if and only if $m$ is odd, there is no loss of information in making this assumption. As a consequence of this convention regarding $m$, we need only consider admissible triples of the form $\{(m,n),(m_j,n_j),(m_j,n_j)\}$, and (again employing (A4) above) this makes it clear that $n$ must be odd as well. To summarize this discussion:
\begin{equation}\label{eq:conventions}
\mbox{\emph{We henceforth adopt the convention that $p,m,n$ are odd integers.}}
\end{equation}
 
For an element $u\in L_{m,n}$ and intertwining operator $\mathcal{Y} \in \mathcal{I}_{(m,n),(m_j,n_j)}^{(m_j,n_j)}$, we express the functions \eqref{eq:tracefcn} as
\begin{align}
 S_{j}(u,\tau) :&= \tr \vert_{L_{m_j,n_j}} o(u) q^{L(0)-\frac{\mathbf{c}}{24}} = q^{h_{m_j,n_j}-\frac{\mathbf{c}}{24}} \sum_{k\geq 0} \tr \vert_{(L_{m_j,n_j})_k} o(u) q^k \label{eq:q1}
\end{align}
with $h_{m_j,n_j}$ and ${\bf c}$ as in \eqr{cw} and \eqr{cc}, respectively. As a special case of \cite[Theorems 4.14 and 5.1]{Miy-Int} (see also \cite{Huang-Int}, \cite[Theorem 5.1]{Yamauchi-Int}), we obtain the following result. 
\begin{thm}\label{thm:linind}
 Let $u$ be a highest weight vector in $L_{m,n}$, and 
 let
 \begin{equation}\label{eq:triples}
 \left\{\{(m,n),(m_j,n_j),(m_j,n_j)\} \vert 1\leq j\leq s\right\}
 \end{equation}
  be the set of all mutually inequivalent admissible triples associated to $(m,n)$. Then the linear space 
\begin{equation}\label{eq:spanSj}
\mc{C}_1^u(L_{m,n})=\mbox{\text{\emph{span}}}_\bb{C}\left\{ S_{j}(u,\tau) \vert 1\leq j\leq s\right\} 
\end{equation}
is invariant under the action \eqr{slash} of the modular group $\G$.\qed
\end{thm} 

We next determine the dimension (over the complex numbers) of this invariant space of $1$-point functions. We note for use below the following result, which follows directly from the definition of admissible triple and the fact that the irreducible $V(p,q)$-modules are uniquely determined by their conformal weight.
\begin{lem}\label{lem:conformal}
 Suppose $\{(m,n),(m_j,n_j),(m_j,n_j)\}$ and $\{(m,n),(m_k,n_k),(m_k,n_k)\}$ are two inequivalent admissible triples. Then $h_{m_j,n_j}\not = h_{m_k,n_k}$.\qed
\end{lem}

Additionally, we will need the following result.
\begin{lem}\label{lem:non0}
 Suppose $\{(m,n),(m_j,n_j),(m_j,n_j)\}$ is an admissible triple and $u$ is a highest weight vector of $L_{m,n}$. Then the leading coefficient of the expansion \eqr{q1} for $S_{j}(u,\tau)$ is nonzero.
\end{lem}

\proof 
 Set $U=L_{m,n}$ and $W=L_{m_j,n_j}$ and denote their $L(0)$-graded eigenspaces by $U=\bigoplus_{k=0}^\infty U_{h_U +k}$ and $W=\bigoplus_{k=0}^\infty W_{h_W +k}$, respectively, where $h_U =h_{m,n}$ and $h_W=h_{m_j,n_j}$. We consider the isomorphism \cite[Corollary 2.13]{Li-DeterminingFusion}
 \[
 \pi \colon \binom{W}{U\, W} \to \text{Hom}_{A(V)}(A(U)\otimes_{A(V)} W_{h_W} ,W_{h_W})
 \]
 where $\pi(\mathcal{Y})(v\otimes b)=o(v)b$ for $v\in A(U)$ and $b\in W_{h_W}$ (see also \cite{FZ-AssociatedTo}). Here $A(U)$ is the $A(V)$-bimodule, where $A(V)$ is the Zhu algebra of $V=V(p,q)$ (see \cite{Li-DeterminingFusion} for a precise definition and more details about $A(V)$, $A(U)$, and the map $\pi$). By our assumption, $N_{(m_i,n_i), (m_j,n_j)}^{(m_k,n_k)}=1$, and we can take $\mathcal{Y} \in \mathcal{I}_{(m_i,n_i), (m_j,n_j)}^{(m_k,n_k)}$ to be nonzero. If the leading coefficient of $S_j (u,\tau)$ is zero then $\tr_{W_{h_W}} o(u)=0$, and since $\dim W_{h_W}=1$, we must have $o(u)=0$ on $W_{h_W}=\C b$ for the highest weight vector $b$. Therefore, $\pi(\mathcal{Y})(u\otimes b)=0$. We claim that $\pi(\mathcal{Y})(v\otimes b)=0$ for all $v\in A(U)$, showing that $\pi (\mathcal{Y})=0$, a contradiction.\\
 \indent We will prove that $o(v)b=0$ for $v\in A(U)$ (by definition $A(U)\subset U$) and $b\in W_{h_W}$ by induction on $\wt v$. Any element in $U_{h_U+k}$ is of the form $v=L(-\ell_1)\cdots L(-\ell_t)u$, where $\ell_1 ,\dots ,\ell_t >0$ and $\ell_1 +\cdots +\ell_t =k$. Additionally, the bracket relations allow us to assume $\ell_j$ is $1$ or $2$ for $1\leq j\leq t$. The case $k=0$ implies $v=u$ and $o(u)b=0$ follows from the discussion above. Suppose $o(v)b=0$ for all $v\in A(U)$ with $\wt v<k$. Consider $v=L(-\ell_1)\cdots L(-\ell_t)u$, where $\ell_1 ,\dots ,\ell_t =1,2$ and $\ell_1 +\cdots +\ell_t =k$ and set $x=L(-\ell_2)\cdots L(-\ell_t)u$ so that $v=L(-\ell_1)x$. Note that $\wt x=h_U +\ell_2 +\cdots +\ell_t$ and $\ell_2 +\cdots +\ell_t<k$, so that by our induction hypothesis $o(x)b=0$. The associator formula gives
 \begin{align*}
 &o\left(L(-\ell_1)x\right)b=\left( \omega (1-\ell_1)x \right)(h_U +\ell_1 +\cdots +\ell_t -1)b \\
 &=\sum_{i\geq 0} (-1)^i \binom{1-\ell_1}{i}\left[ \omega (1-\ell_1 -i)x(\wt x-1+\ell_1 +i) -(-1)^{1-\ell_1}x(\wt x-i)\omega (i) \right]b\\
 &=(-1)^{\ell_1} \left[\binom{1-\ell_1}{0}x(\wt x)L(-1)b+h_W\binom{1-\ell_1}{1}x(\wt x-1)b\right] =0
 \end{align*}
 since
 \begin{align*}
 &x(\wt x)L(-1)b=L(-1)x(\wt x)b-[\omega (0),x(\wt x)]b=-(L(-1)x)(\wt x)b =\wt x o(x)b=0.
 \end{align*}
 This completes the proof of the lemma. \qed
 
 \begin{rmk}
 It is an immediate corollary that under the assumptions of the previous lemma, $S_j (u,\tau) \not =0$ for a highest weight element $u\in L_{m,n}$ (see also \cite[Lemma 6.2]{Miy-Int}). 
 \end{rmk}
 
\begin{cor}\label{cor:linInd}
 Let $u$ be a highest weight vector for $L_{m,n}$ and suppose
 \[
 \left\{\{(m,n),(m_j,n_j),(m_j,n_j)\} \vert 1\leq j\leq s\right\}
 \]
 is the associated set of inequivalent admissible triples. Then the $1$-point functions 
  \[
  \left\{S_{j}(u,\tau) \vert 1\leq j\leq s\right\} \]
are linearly independent.
\end{cor}

\proof 
 By Lemma \ref{lem:conformal} it suffices to show that the leading coefficient in the $q$-expansions of any $S_{j}(u,\tau)$ is nonzero, and this fact is provided by Lemma \ref{lem:non0}. \qed
 
\medskip
This corollary implies that the dimension of the space \eqr{spanSj} is equal to the integer $s$ in Theorem \ref{thm:linind}, and we now give a formula for this integer.
\begin{lem}\label{lem:rep}
Given $m$ and $n$, the set of integer pairs $(m_j,n_j)$ satisfying 
\begin{equation}\label{eq:bounds} 
\frac{p+1}{2}\leq m_j\leq p-\frac{m+1}{2} \quad \mbox{and} \quad \frac{n+1}{2}\leq n_j\leq q-\frac{n+1}{2}
\end{equation}
gives a complete set of admissible triples of the form
$\{(m,n),(m_j,n_j),(m_j,n_j)\}$.
\end{lem}

\proof Conditions (A2) and (A3) above show that the inequalities 
\[
\frac{m}{2}<m_j<p-\frac{m}{2} \quad \mbox{and} \quad \frac{n}{2}<n_j<q-\frac{n}{2}
\]
hold, and $p,m,n$ are odd by convention \eqr{conventions}. This implies the inequalities
\[
\frac{m+1}{2}\leq m_j\leq p-\frac{m+1}{2} \quad \mbox{and} \quad \frac{n+1}{2}\leq n_j\leq q-\frac{n+1}{2}.
\]
The final modification follows from observing that $\frac{p-1}{2}<m_j$ if and only if $p-m_j\geq\frac{p+1}{2}$. Since the modules $L_{m_j,n_j}$ and $L_{p-m_j,q-n_j}$ coincide and (A1) implies $\frac{m+1}{2}<\frac{p+1}{2}$, one way to ensure that we only count once a given $L_{m_j,n_j}$ acted on by $L_{m,n}$ is to change the lower bound on $m_j$ to $\frac{p+1}{2}$ and let $n_j$ vary as indicated above. This completes the proof.\qed

\medskip
Again keeping in mind our conventions \eqr{conventions}, we now obtain the desired dimension formula.
\begin{cor}\label{cor:dim}
Given an irreducible $V(p,q)$-module $L_{m,n}$ and a highest weight vector $u\in L_{m,n}$, the dimension of the associated space \eqr{spanSj} of $1$-point functions is $\frac{(p-m)(q-n)}{2}$.
\end{cor} 

\proof From \eqr{bounds} we see that the number of possible $m_j$ is
\[
p-\frac{m+1}{2}-\frac{p+1}{2}+1=\frac{p-m}{2}
\]
and the number of possible $n_j$ is 
\[
q-\frac{n+1}{2}+1-\frac{n+1}{2}=q-n.
\]
Since the $m_j$ and $n_j$ are completely independent, we multiply these numbers together to obtain all possible pairs $(m_j,n_j)$ appearing in admissible triples of the form given in Lemma \ref{lem:rep}. Since each corresponding fusion coefficient $N_{(m,n)\ (m_j,n_j)}^{(m_j,n_j)}$ is $1$, by Lemma \ref{lem:non0} there is a unique nonzero $1$-point function $S_j(u,\tau)$ associated to each pair $(m_j,n_j)$. These form a linearly independent set by Corollary \ref{cor:linInd}, so the dimension of \eqr{spanSj} is equal to the number of such pairs.\qed 

\medskip
Since each irreducible $V(p,q)$-module $L_{m,n}$ is generated by the Virasoro modes acting on its highest weight vector $u$, we obtain the following important result.
 
\begin{thm}\label{thm:1ptgen}
Let $u$ be a highest weight vector of $L_{m,n}$. Then the space of $1$-point functions $\mathcal{C}_1 (L_{m,n})$ associated to $L_{m,n}$ is equal to the space $\mathcal{R} \mathcal{C}_1^u (L_{m,n})$, where $\mc{C}_1^u(L_{m,n})$ is as in \eqr{spanSj} and $\mc{R}$ denotes the ring \eqr{R} of modular differential operators.
\end{thm}

\proof Set $W=L_{m,n}$ and let the notation be as in Corollary \ref{cor:linInd}. Consider a function $f(\tau) \in \mathcal{C}_1(W)$. Then $f(\tau)$ is a linear combination of functions of the form $S_j(v_j,\tau)$ for $v_j\in W$ and $1\leq j\leq s$. Meanwhile, each $v_j$ can be written as a linear combination of homogeneous elements $w_r$ with respect to the $L[0]$-grading, say $w_r \in W_{[d_r]}$ where $d_r \in \mathbb{N}$. Therefore it suffices to show $S_j(w_r,\tau)\in \mathcal{R} \mathcal{C}_1^u(W)$. Since $w_r \in W_{[d_r]}$, there exists some $\ell_r \in \mathbb{N}$ and $n_t=1,2$ for $1\leq t\leq \ell_r$, such that $w_r =L[-n_1] \cdots L[-n_{\ell_r}]u$. By Lemma \ref{lem:surjective1} there exists $m_j \in \mathbb{N}$ and $\phi_t \in \mathcal{M}$ for $0\leq t\leq m_r$, such that $S_j(L[-n_1] \cdots L[-n_{\ell_r}]u, \tau) = \left(\phi_0 +\phi_1 \partial +\cdots +\phi_{m_r} \partial^{m_r}\right)S_j(u,\tau)$. Thus, 
 \begin{equation}\label{eq:Rsum}
 S_j(w_r ,\tau)=\left(\phi_0 +\phi_1 \partial +\cdots +\phi_{m_r} \partial^{m_r}\right)S_j(u,\tau).
 \end{equation}
  Moreover, following the proof of Proposition 2(b) in \cite{DM-Higher} we find the same $m_r \in \mathbb{N}$, $\phi_t \in \mathcal{M}$ for $0\leq t\leq m_r$, and equality in \eqref{eq:Rsum} hold for every $1\leq j\leq s$. Therefore, $S_j(w_r ,\tau)=\left(\phi_0 +\phi_1 \partial +\cdots +\phi_{m_r} \partial^{m_r}\right)S_j(u,\tau) \in \mathcal{R}\mathcal{C}_1^u(W)$, and we conclude $f(\tau) \in \mathcal{R}\mathcal{C}_1^u(W)$.
 \qed 
 
\medskip
This result has the following practical consequence. Given an irreducible $V(p,q)$-module $L_{m,n}$ with highest weight vector $u$, together with the set \eqr{triples} of inequivalent admissible triples, Theorem \ref{thm:linind} tells us there is a matrix representation
\begin{equation}\label{eq:rhomn}
\rho_{m,n}:\G\rightarrow\gln{s}
\end{equation}
such that
\begin{equation}\label{eq:virF}
F(\tau)=\cvec{S_1(u,\tau)}{\vdots}{S_s(u,\tau)}=\cvec{q^{\lambda_1} \sum_{k\geq 0} \tr \vert_{(L_{m_1,n_1})_k} o(u) q^k}{\vdots}{q^{\lambda_s} \sum_{k\geq 0} \tr \vert_{(L_{m_s,n_s})_k} o(u) q^k}
\end{equation}
is a weakly holomorphic vector-valued modular form of weight $\wt[u]$ for $\rho_{m,n}$. Here (and in what follows) we set
\begin{equation}\label{eq:lambdaj}
\lambda_j:=h_{m_j,n_j}-\frac{\mathbf{c}}{24}
\end{equation}
and recall that $\bf{c}$ and $h_{m_j,n_j}$ are given by \eqr{cc} and \eqr{cw}, respectively, $s=\dim\rho_{m,n}=\frac{(p-m)(q-n)}{2}$ by Corollary \ref{cor:dim}, and our choice of multiplier system \eqr{upT} (in weight $k=\wt[u]$) is in effect. Evidently we have
\begin{equation}\label{eq:rhoTmn}
\rho_{m,n}(T)=\diag{\e{r_1},\dots,\e{r_s}}
\end{equation}
with $r_j=h_{m_j,n_j}-\frac{\mathbf{c}}{24}-\frac{h_{m,n}}{12}$.
Letting the ring \eqr{R} of modular differential operators act componentwise on $F$, Theorem \ref{thm:1ptgen} and the covariance \eqr{covar} of the modular derivative imply that the space $\mc{C}_1(L_{m,n})$ is organized into a $\bb{Z}$-graded space 
\begin{equation}\label{eq:Vrho}
\mc{V}(\rho_{m,n})=\mc{R}F=\bigoplus_{k\geq0}\mc{V}(\wt[u]+2k,\rho_{m,n})
\end{equation} 
of weakly holomorphic vector-valued modular forms for $\rho_{m,n}$, such that for any $f\in\mc{C}_1(L_{m,n})$ there is a vector $G\in\mc{V}(\rho_{m,n})$ whose components contain $f$ in their span. In the subsections below, we compare $\mc{V}(\rho_{m,n})$ to the space of holomorphic vector-valued modular forms for $\rho_{m,n}$. In order to facilitate this we utilize results \cite{Marks1, MarksMason} of Mason and the second author, thus we review said results before proceeding to the proofs. 

Suppose that $L_{m,n}$ is an irreducible $V(p,q)$-module such that the associated representation \eqr{rhomn} is irreducible of dimension $s\leq3$, with $\rho_{m,n}(T)$ diagonal as in \eqr{rhoTmn}. Let $u$ be a highest weight vector for $L_{m,n}$, set $k=\wt[u]$, and let $\up$ be the multiplier system \eqr{upT}. For $1\leq j\leq s$ let $\alpha_j$ denote the unique real number such that 
\begin{equation}\label{eq:alphaj}
0\leq\alpha_j<1 \quad \mbox{and} \quad \alpha_j\equiv r_j+\frac{k}{12}\pmod{\bb{Z}}.
\end{equation} 
Then the space of holomorphic vector-valued modular forms for $\rho_{m,n}$ and $\up$ has a $2\bb{Z}$-grading
\begin{equation}\label{eq:Hrhomn}
\mc{H}(\rho_{m,n})=\bigoplus_{l\geq0}\mc{H}_{k_0+2l}(\rho_{m,n}),
\end{equation}
where 
\begin{equation}\label{eq:k0}
k_0=\frac{12\sum_j\alpha_j}{s}+1-s
\end{equation}
gives the minimal weight of any nonzero holomorphic vector-valued modular form for $\rho_{m,n}$ and $\up$. Furthermore, there is a vector 
\begin{equation}\label{eq:F0}
F_0=\cvec{f_1}{\vdots}{f_s}=\cvec{q^{\alpha_1}+\sum_{n\geq0}a_1(n)q^{\alpha_1+n}}{\vdots}{q^{\alpha_s}+\sum_{n\geq0}a_s(n)q^{\alpha_s+n}}
\end{equation}
such that $\mc{H}_{k_0}(\rho_{m,n})=\bb{C}F_0$ and $\mc{H}(\rho_{m,n})=\mc{R}F_0$ is the cyclic $\mc{R}$-module generated by $F_0$, where $\mc{R}$ denotes the ring \eqr{R} of modular differential operators. Theorem \ref{thm:1ptgen} shows that $\mc{V}(\rho_{m,n})$, the associated space \eqr{Vrho} of vector-valued $1$-point functions for $L_{m,n}$ is also a cyclic $\mc{R}$-module, and from this it follows that $\mc{V}(\rho_{m,n})$ contains $\mc{H}(\rho_{m,n})$ if and only if the minimal weight vector $F_0$ for $\mc{H}(\rho_{m,n})$ is contained in $\mc{V}(\rho_{m,n})$. A special case of this occurs when $F_0$ is a scalar multiple of the generator \eqr{virF} for $\mc{V}(\rho_{m,n})$, so that $\mc{V}(\rho_{m,n})$ and $\mc{H}(\rho_{m,n})$ coincide. In order for this to occur, it is necessary and sufficient that the $\alpha_j$ in \eqr{F0} are equal to the leading exponents of \eqr{virF} and that the weight $k_0$ in \eqr{k0} is equal to $k=\wt[u]$. In any event, if all the exponents of \eqr{virF} are nonnegative then \eqr{virF} is holomorphic and $\mc{V}(\rho_{m,n})$ is contained in $\mc{H}(\rho_{m,n})$, and this containment is proper if and only if at least one of these exponents is greater than or equal to one. 

We note here, for use in what follows, that the leading exponents of \eqr{virF} are given by
\begin{eqnarray}\label{eq:lambdas}
\lambda_j=h_{m_j,n_j}-\frac{\mathbf{c}}{24}&=&\frac{(n_jp-m_jq)^2-(p-q)^2}{4pq}-\frac{1}{24}\left[1-\frac{6(p-q)^2}{pq}\right] \notag \\
&=&\frac{(n_jp-m_jq)^2}{4pq}-\frac{1}{24}. 
\end{eqnarray}
\begin{prop}\label{prop:primedim}
Suppose that $s=\dim\rho_{m,n}$ is equal to 1 or a prime number. Then one of the following holds:
\begin{enumerate}
\item $m=p-2$, $n=q-s$, and 
\begin{equation}\label{eq:case1lambdaj}
\lambda_j=\frac{3(1+s-2j)^2p^2+2(2+3s-6j)pq+3q^2}{48pq},\ 1\leq j\leq s,
\end{equation}
whereas the exponents of \eqr{rhoTmn} are given by 
\begin{equation}\label{eq:case1rj}
r_j=\frac{(3[2j-(s+1)]^2+1-s^2)p+2(1+5s-6j)q}{48q},\ 1\leq j\leq s.
\end{equation}

\item $m=p-2s$, $n=q-1$, and 
\begin{equation}\label{eq:lambdajcase2}
\lambda_j=\frac{3(1-2j)^2q-2p}{48p},\ 1\leq j\leq s,
\end{equation}
whereas the exponents of \eqr{rhoTmn} are given by 
\begin{equation}\label{eq:case2rj}
r_j=\frac{[3(1-2j)^2+1-4s^2]q+4(s-1)p}{48p},\ 1\leq j\leq s.
\end{equation}
\end{enumerate}
These cases coincide exactly when $s=1$. In each case, the  formula
\begin{equation}\label{eq:monic}
h_{m,n}=\frac{12\sum_{j=1}^s\lambda_j}{s}+1-s
\end{equation}
obtains. 
\end{prop}

\proof From Corollary \ref{cor:dim} we have two cases (which coincide for $s=1$), namely $\frac{p-m}{2}=1$, $q-n=s$ or $\frac{p-m}{2}=s$, $q-n=1$.

Considering the former case first, we have $m=p-2$, $n=q-s$. From Lemma \ref{lem:rep} we obtain $m_j=\frac{p+1}{2}$ and $n_j=\frac{n+2j-1}{2}$ for $1\leq j\leq s$. Substituting this into \eqr{lambdas} gives \eqr{case1lambdaj}. In this case we have
\begin{equation}\label{eq:hmn}
\frac{h_{m,n}}{12}=\frac{(np-mq)^2-(p-q)^2}{48pq}=\frac{3q^2+2(1-2s)pq+(s^2-1)p^2}{48pq},
\end{equation}
and setting $r_j=\lambda_j-\frac{h_{m,n}}{12}$ yields \eqr{case1rj}. We next compute the right side of \eqr{monic} and obtain
\begin{eqnarray*}
1-s+\frac{12}{s}\sum_{j=1}^s\lambda_j&=&1-s+\frac{1}{4spq}\left[3p^2\sum_{j=1}^s[(1+s)^2-4(1+s)j+4j^2]\right.\\
&\ &+\left.2pq\sum_{j=1}^s(2+3s-6j)+3sq^2\right]\\
&=&1-s+\frac{1}{4spq}\left[3p^2\left(s(1+s)^2-2s(1+s)^2+\frac{2s(s+1)(2s+1)}{3}\right)\right.\\
&\ &+\left.2pq(s(2+3s)-3s(s+1))+3sq^2\right]\\
&=&1-s+\frac{1}{4pq}[(s^2-1)p^2-2pq+3q^2]\\
&=&\frac{(s^2-1)p^2+2(1-2s)pq+3q^2}{4pq}.
\end{eqnarray*}
Comparing this with \eqr{hmn} shows that \eqr{monic} holds.

In the second case, we have $m=p-2s$, $n=q-1$, so from \eqr{bounds} we have $m_j=\frac{p+2j-1}{2}$, $n_j=\frac{n+1}{2}=\frac{q}{2}$ for $1\leq j\leq s$. Using this in \eqr{lambdas} now yields \eqr{lambdajcase2}. In this setting we have
\begin{equation}\label{eq:hmn2}
\frac{h_{m,n}}{12}=\frac{[(q-1)p-(p-2s)q]^2-(p-q)^2}{48pq}=\frac{(4s^2-1)q+2(1-2s)p}{48p},
\end{equation}
whereas the right side of \eqr{monic} yields in this case
\begin{eqnarray*}
1-s+\frac{12}{s}\sum_{j=1}^s\lambda_j&=&1-s+\frac{1}{4sp}\sum_{j=1}^s\left(3(1-2j)^2q-2p\right)\\
&=&1-s+\frac{1}{4sp}\left[3\left(s-2s(s+1)+\frac{2s(s+1)(2s+1)}{3}\right)q-2sp\right]\\
&=&1-s+\frac{(4s^2-1)q-2p}{4p}\\
&=&\frac{(4s^2-1)q+2(1-2s)p}{4p}.
\end{eqnarray*}
Comparing this to \eqr{hmn2} shows that \eqr{monic} holds in this case as well.\qed 

\medskip
To summarize, it is now apparent from \eqr{k0} and \eqr{monic}, together with the discussion preceding Proposition \ref{prop:primedim}, that the following result holds.
\begin{cor}\label{cor:equal}
Let $\lambda_j=h_{m_j,n_j}-\frac{\mathbf{c}}{24}$ as in Proposition \ref{prop:primedim}, and suppose that $\rho_{m,n}$ in \eqr{rhomn} is irreducible of dimension less than four. Then the space $\mc{V}(\rho_{m,n})$ of vector-valued $1$-point functions is equal to the space $\mc{H}(\rho_{m,n})$ of holomorphic vector-valued modular forms if and only if $0\leq\lambda_j<1$ for $1\leq j\leq s$. If $\lambda_j\geq0$ for each $j$ and $\lambda_j\geq1$ for at least one $j$ then $\mc{V}(\rho_{m,n})$ is properly contained in $\mc{H}(\rho_{m,n})$.\qed
\end{cor}

We are now properly situated to analyze the low-dimension setting of modular group representations and spaces of $1$-point functions arising from Virasoro minimal model intertwining operators, and proceed on a dimension-by-dimension basis as follows.


\subsection{Dimension one}

If the dimension of $\rho_{m,n}$ is one, then the two cases of Proposition \ref{prop:primedim} coincide and we have the VOA  $V=V(p,q)$ with $p\geq3$ odd and $q\geq2$ even. From Corollary \ref{cor:dim} one sees that the relevant $V$-module in this setting is  $L_{m,n}=L_{p-2,q-1}$, and both \eqr{case1rj} and \eqr{case2rj} imply that $\rho(T)=1$. Thus $\rho$ is  the trivial one-dimensional representation of $\G$ (regardless of $p$ and $q$) and the associated space $\mc{V}(\rho_{p-2,q-1})$ of $1$-point functions consists of weakly holomorphic modular forms for $\G$. The generator \eqr{virF} for $\mc{V}(\rho_{p-2,q-1})$ has leading exponent 
$$\lambda_1=\frac{3q-2p}{48p},$$
and we note that $0\leq\lambda_1<1$ if and only if 
\begin{equation}\label{eq:dim1}
\frac{2}{3}\leq\frac{q}{p}<\frac{50}{3}.
\end{equation}
By Corollary \ref{cor:equal}, we obtain the following result.
\begin{thm}\label{thm:dim1}
Suppose $L_{m,n}$ is an irreducible $V(p,q)$-module such that the modular group representation $\rho_{m,n}$ associated to $L_{m,n}$ is one-dimensional. Then (recalling the conventions \eqr{conventions}) $p\geq3$ is odd, $q\geq2$ is even, $(m,n)=(p-2,q-1)$, and the space $\mc{V}(\rho_{m,n})$ of vector-valued 1-point functions associated to $L_{m,n}$ is equal to the space $\mc{H}(\rho_{m,n})$ of holomorphic vector-valued modular forms for $\rho_{m,n}$ exactly when $\frac{q}{p}$ satisfies \eqr{dim1}.\qed
\end{thm}

We recall here that the kernel of any one-dimensional representation of $\G$ is necessarily congruence, since the commutator subgroup $\G'$ of $\G$ is congruence of level 12 (see e.g.\ \cite[Theorem 1.3.1]{Rankin}). In fact $\G/\G'$ and $\Hom(\G,\cstar)$ are each cyclic of order 12, with the latter group being generated by $\chi:\G\rightarrow\cstar$ where $\chi(T)=\e{\frac{1}{12}}$. Since the determinant of any representation of $\G$ defines an element of $\Hom(\G,\cstar)$, it follows a representation $\rho$ of $\G$ is irreducible if no subproduct of the eigenvalues of $\rho(T)$ is a twelfth root of unity. We will make use of this fact in the following subsections.


\subsection{Dimension two}
Suppose now that $L_{m,n}$ is an irreducible $V=V(p,q)$-module such that \eqr{rhomn} is two-dimensional. In the first case of Proposition \ref{prop:primedim} (and keeping in mind our conventions \eqr{conventions}) we find that $p\geq3$ and $q\geq3$ are both odd, and the $V$-module $L_{p-2,q-2}$ yields a two-dimensional representation $\rho=\rho_{p-2,q-2}$ of $\G$ such that the exponents of \eqr{rhoTmn} are
$$r_1=\frac{5}{24},\ r_2=-\frac{1}{24}.$$
Since neither eigenvalue of \eqr{rhoTmn} is a twelfth root of unity, we see that in all cases $\rho_{p-2,q-2}$ is the same irreducible representation, regardless of $p$ and $q$, and is in fact \cite[Table 3]{Mason} congruence. The leading exponents of \eqr{virF} are given by 
$$\lambda_1=\frac{3p^2+4pq+3q^2}{48pq},\ \lambda_2=\frac{3p^2-8pq+3q^2}{48pq},$$
and one may check that $0\leq\lambda_1,\lambda_2<1$ if and only if 
\begin{equation}\label{eq:dim2case1}
\frac{22-5\sqrt{19}}{3} \leq \frac{q}{p} < \frac{4-\sqrt{7}}{3} \quad \mbox{or} \quad \frac{4+\sqrt{7}}{3}\leq\frac{q}{p}<\frac{22+5\sqrt{19}}{3}.
\end{equation}

In the second case of Proposition \ref{prop:primedim}, we have  $p\geq5$ odd and $q\geq2$ even, and the $V$-module $L_{p-4,q-1}$ yields a two-dimensional representation $\rho=\rho_{p-4,q-1}$ of $\G$ whose exponents in \eqr{rhoTmn} are 
$$r_1=\frac{p-3q}{12p},\ r_2=\frac{p+3q}{12p}.$$
It is evident that for all $p$ and $q$ under consideration $\rho$ is irreducible of level $12p$, since once again neither eigenvalue of \eqr{rhoTmn} is a twelfth root of unity. From this observation and \cite{Mason}, it follows that for all $p>5$ $\rho$ is infinite image (thus noncongruence), whereas for $p=5$ we obtain (by varying $q$) four of the level 60 congruence representations tabulated in \cite{Mason}.

The leading exponents of \eqr{F} in this case are
$$\lambda_1=\frac{3q-2p}{48p},\ \lambda_2=\frac{27q-2p}{48p},$$
and from this it follows that $0\leq\lambda_1,\lambda_2<$ if and only if
\begin{equation}\label{eq:dim2case2}
\frac{2}{3}\leq\frac{q}{p}<\frac{50}{27}.
\end{equation}
We summarize the above findings in the following theorem.
\begin{thm}\label{thm:dim2}
Suppose that $L_{m,n}$ is an irreducible $V(p,q)$-module such that the associated representation $\rho_{m,n}$ of $\G$ is two-dimensional. Then $\rho_{m,n}$ is irreducible and (keeping in mind the conventions \eqr{conventions}) one of the following two cases obtains:
\begin{enumerate}
\item[(i)]$p,q\geq5$ are odd and $(m,n)=(p-2,q-2)$.
\item[(ii)]$p\geq5$ is odd, $q\geq2$ is even, and $(m,n)=(p-4,q-1)$.
\end{enumerate}
In case $(i)$ $\rho_{m,n}$ is the same congruence representation regardless of $p$ and $q$, with $\rho_{m,n}(T)=\diag{\e{\frac{5}{24}},\e{-\frac{1}{24}}}$, and the space $\mc{V}(\rho_{m,n})$ of vector-valued 1-point functions associated to $L_{m,n}$ is equal to the space $\mc{H}(\rho_{m,n})$ of holomorphic vector-valued modular forms for $\rho_{m,n}$ exactly when $\frac{q}{p}$ satisfies \eqr{dim2case1}. In case $(ii)$ $\rho_{m,n}$ is congruence if and only if $p=5$, and $\mc{V}(\rho_{m,n})=\mc{H}(\rho_{m,n})$ exactly when $\frac{q}{p}$ satisfies \eqr{dim2case2}. 
\end{thm}


\subsection{Dimension three}

Suppose $V=V(p,q)$ with $p\geq3$ odd and $q\geq4$ even. Then the first case of Proposition \ref{prop:primedim} holds and the action of the $V$-module $L_{p-2,q-3}$ yields a three-dimensional representation $\rho=\rho_{p-2,q-3}$ of the modular group. From \eqr{case1rj} we obtain
$$r_1=\frac{p+5q}{12q}\ ,r_2=\frac{q-p}{6q},\ r_3=\frac{p-q}{12q}.$$
From this it follows that $\rho$ is irreducible (for all $p$ and $q$ under consideration) of level $12q$. Since $r_1-r_3=\frac{1}{2}$, it follows from \cite[Proposition 5.1]{Marks2} that each of these representations is finite image. Furthermore, it follows from \cite[Corollary 3.5]{Marks2} that these representations have a noncongruence kernel if $12q$ does not divide $2^8\cdot3^4\cdot5^2\cdot7^2$. From \eqr{case1lambdaj} we obtain
$$\lambda_1=\frac{12p^2+10pq+3q^2}{48pq},\ \lambda_2=\frac{3q-2p}{48p},\ \lambda_3=\frac{12p^2-14pq+3q^2}{48pq},$$
and from this it follows that $0\leq\lambda_1,\lambda_2,\lambda_3<1$ if and only if
\begin{equation}\label{eq:dim3}
\frac{2}{3} \leq \frac{q}{p} < \frac{7-\sqrt{13}}{3} \quad \mbox{or} \quad
 \frac{7+\sqrt{13}}{3} \leq \frac{q}{p} < \frac{19+5\sqrt{13}}{3}.
\end{equation}

Now assume that $V=V(p,q)$ with $p\geq7$ odd and $q\geq2$ even. Then the second case of Proposition
\ref{prop:primedim} holds and the $V$-module $L_{p-6,q-1}$ yields a three-dimensional representation $\rho=\rho_{p-6,q-1}$ of the modular group. We have from \eqr{case2rj} that
$$r_1=\frac{p-4q}{6p},\ r_2=\frac{p-q}{6p},\ r_3=\frac{p+5q}{6p}.$$
This shows that (for any $p$ and $q$ in the range of consideration) $\rho$ is irreducible of level $6p$, and from \cite[Proposition 5.1]{Marks2} it follows that $\rho$ is infinite image (thus noncongruence). From \eqr{lambdajcase2} we obtain
$$\lambda_1=\frac{3q-2p}{48p},\ \lambda_2=\frac{27q-2p}{48p},\ \lambda_3=\frac{75q-2p}{48p},$$
there are no values of $p$ and $q$ such that $0\leq\lambda_j<1$ for each $j$. Indeed, for this to hold for $\lambda_1$ it would be necessary that $\frac{q}{p}\geq\frac{2}{3}$, while for $\lambda_3$ to satisfy this condition the inequality $\frac{q}{p}<\frac{2}{3}$ must be satisfied. Thus in this case it will never happen that $\mc{V}(\rho_{p-6,q-1})=\mc{H}(\rho_{p-6,q-1})$.
We summarize our findings in the following result.
\begin{thm}\label{thm:dim3}
Suppose that $L_{m,n}$ is an irreducible $V(p,q)$-module such that the associated representation $\rho_{m,n}$ of $\G$ is three-dimensional. Then $\rho_{m,n}$ is irreducible and (keeping in mind the conventions \eqr{conventions}) one of the following two cases obtains:
\begin{enumerate}
\item[(i)]$p\geq3$ is odd, $q\geq4$ is even, $(m,n)=(p-2,q-3)$.
\item[(ii)]$p\geq7$ is odd, $q\geq2$ is even, and $(m,n)=(p-6,q-1)$.
\end{enumerate}
In case $(i)$ $\rho_{m,n}$ is finite image of level $12q$, $\ker\rho_{m,n}$ is noncongruence if $q$ does not divide $2^6\cdot3^3\cdot5^2\cdot7^2$, and the space $\mc{V}(\rho_{m,n})$ of vector-valued 1-point functions associated to $L_{m,n}$ is equal to the space $\mc{H}(\rho_{m,n})$ of holomorphic vector-valued modular forms for $\rho_{m,n}$ exactly when $\frac{q}{p}$ satisfies \eqr{dim3}. In case $(ii)$ $\rho_{m,n}$ is infinite image (thus noncongruence), and $\mc{V}(\rho_{m,n})$ does not coincide with $\mc{H}(\rho_{m,n})$ for any $p$ and $q$.\qed
\end{thm}

\subsection{Arbitrary dimension}\label{subsec:arb}

Again we consider an irreducible $V(p,q)$-module $L_{m,n}$ with highest weight vector $u$, and the corresponding vector \eqr{virF} of 1-point functions. Setting $k=\wt[u]=h_{m,n}$ in \eqr{upT} and using  \eqr{cw} gives 
\[
\up(T)=\e{\frac{h_{m,n}}{12}}=\e{\frac{(np-mq)^2-(p-q)^2}{48pq}}.
\]
Again employing \eqr{cw} and also \eqr{cc}, the $j^{\text{th}}$ exponent of \eqr{rhoTmn} reads 
\begin{equation}\label{eq:rj}
r_j=h_{m_j,n_j}-\frac{\mathbf{c}}{24}-\frac{h_{m,n}}{12}=\frac{12(n_jp-m_jq)^2-(np-mq)^2+(p-q)^2-2pq}{48pq}.
\end{equation}
Writing $r_j=\frac{x_j}{48pq}$, the identities
\begin{equation}\label{eq:pq}
x_j\equiv12n_j^2-n^2+1\pmod{q} \quad \mbox{and} \quad  x_j\equiv12m_j^2-m^2+1\pmod{p}
\end{equation}
obtain, which will be utilized in what follows. 

We will now give arithmetic conditions on $m$ and $n$ that imply the representation $\rho=\rho_{m,n}$ in \eqr{rhomn} is noncongruence. The basic idea is that one may deduce that $\rho$ is noncongruence by comparing $\dim\rho$ to the \emph{level} of $\rho$, which by definition \cite{W} is the order $N$ of $\rho(T)$ in $\rho(\G)$. In fact, for each $d\geq 1$ there is a finite set of primes $S=S(d)$ such that any prime dividing the level of a $d$-dimensional congruence representation of $\G$ is an element of $S$. Put another way, for $\rho$ to be noncongruence it is sufficient (though generally not necessary) that there exists a prime divisor of $N$ that is not in $S$. This is the strategy we employ below. That is, by comparing the level and the dimension of certain $\rho$ we will be able to show that $\ker\rho$ is a noncongruence subgroup of $\G$.

To facilitate this, we recall that \eqr{gammaN}, the principal congruence subgroup of level $N$, is by definition the kernel of the reduction$\pmod{N}$ map $\G\rightarrow \mbox{SL}_2(\bb{Z}/N\bb{Z})$. Thus, if $\rho:\G\rightarrow\gln{d}$ is congruence of level $N$, then $\rho(\G)$ is isomorphic to a quotient of $\mbox{SL}_2(\bb{Z}/N\bb{Z})$. Since $\mbox{SL}_2(\bb{Z}/N\bb{Z})\cong\prod \mbox{SL}_2(\bb{Z}/r_i^{t_i}\bb{Z})$, where 
\begin{equation}\label{eq:Ndecomp}
N=\prod r_i^{t_i}
\end{equation} 
is the prime decomposition of $N$, it suffices to know the minimal dimensions of irreducible representations of $\mbox{SL}_2(\bb{Z}/r^t \bb{Z})$ for any prime $r$ and integer $t\geq 1$. Fortunately this is well-known, and we record here some results of Nobs and Nobs-Wolfart, which are further elucidated (in English) in the doctoral thesis of Eholzer \cite{E}.
\begin{thm}[\cite{Nobs, Nobs-Wolfart}]\label{thm:NW}
Let $r$ be prime and $t$ a positive integer. Then the smallest dimension of a nontrivial irreducible representation of $\mbox{SL}_2\left(\bb{Z}/r^t \bb{Z}\right)$ is
\begin{enumerate}
\item[(i)] $1$ for $t \le 2$, $r=2$,
\item[(ii)] $2$ for $t=3$, $r=2$,
\item[(iii)] $2^{t-4}\cdot3$ for $t\geq4$, $r=2$,
\item[(iv)] $\frac{r-1}{2}$ for $t=1$, $r>2$, and
\item[(v)] $\frac{(r^2-1)r^{t-2}}{2}$ for $t>1$, $r>2$. \qed
\end{enumerate}
\end{thm}

If $N$ is composite as above, then any representation $\rho$ of $\mbox{SL}_2(\bb{Z}/N\bb{Z})$ will be isomorphic to a tensor product $\otimes\rho_i$, where $\rho_i$ is a representation of $\mbox{SL}_2(\bb{Z}/r_i^{t_i}\bb{Z})$. Thus the minimal dimension of a congruence $\rho$ of level $N$ will be given by the product of the minimal possible dimensions of the $\rho_i$, as given by Theorem \ref{thm:NW}.

For a prime number $r$, we will employ the standard notation $\nu_r(x)$ to denote the $r$-adic valuation of the integer $x\neq0$, and the extension $\nu_r\left(\frac{a}{b}\right)=\nu_r(a)-\nu_r(b)$ to any nonzero $\frac{a}{b}\in\bb{Q}$. In all that follows below, we assume that $L_{m,n}$ is an irreducible module for a minimal model $V(p,q)$, and $\rho=\rho_{m,n}$ is the representation \eqr{rhomn} with level $N$ given by \eqr{Ndecomp}. We denote by $(m_j,n_j)$, $1\leq j\leq s=\dim\rho_{m,n}$ the integer pairs appearing in \eqr{pq}.

\begin{lem}\label{lem:pdiv}
Suppose $r>3$ is a prime dividing $p$ and $m\leq p-4$. Then $\nu_r(N)=\nu_r(p)$.
\end{lem}

\proof Since $p-m\geq4$, both $m_1:=\frac{p+1}{2}$ and $m_2:=m_1+1$ satisfy \eqr{bounds}. It then follows from \eqr{pq} that, regardless of $n_1$, r divides the numerator $x_1$ of the exponent $r_1$ in \eqr{rj} if and only if $r$ divides $12m_1^2-m^2+1$. If $r$ does not divide this quantity, then $r$ divides the level of $\rho_{m,n}$, and $\nu_r(N)=\nu_r(p)$. Assume that $r$ does divide $x_1$. Again it follows from \eqr{pq} that $r$ divides the numerator of $r_2$ if and only if $r$ divides
$$12m_2^2-m^2+1=12(m_1+1)^2-m^2+1=(12m_1^2-m^2+1)+12(2m_1+1).$$
Thus if $r$ divides the numerator of $r_2$ then $r$ divides $12(2m_1+1)$. Since $r>3$ this would imply that $r$ divides $2m_1+1=p+2$. But this would imply that $r$ divides 2, an impossibility. Thus in this case we also have that $r$ divides $N$, and $\nu_r(N)=\nu_r(p)$.\qed

\medskip
Similarly, we have the following result.
\begin{lem}\label{lem:qdiv}
Suppose $r>3$ is a prime dividing $q$ with $n\leq q-3$. Then $\nu_r(N)=\nu_r(q)$.
\end{lem}

\proof  Since $q-n\geq3$, $n_1:=\frac{n+1}{2}$, $n_2:=n_1+1$, and $n_3:=n_1+2$ each satisfy \eqr{bounds}. If $r$ does not divide $12n_1^2-n^2+1$, then it follows from \eqr{pq} that $r$ does not divide the numerator of $r_1$ in \eqr{rj}, and therefore $r$ divides the level $N$ of $\rho_{m,n}$ with $\nu_r(N)=\nu_r(q)$. Assume that $r$ does divide the numerator of $r_1$. Then from \eqr{pq} we see that $r$ divides the numerator of $r_2$ if and only if $r$ divides
\[
12n_2^2-n^2+1=12(n_1+1)^2-n^2+1=(12n_1^2-n^2+1)+12(2n_1+1).
\]
Since $r>3$, this shows that $r$ must divide $2n_1+1=n+2$ if $r$ divides the numerator of $r_2$. If this does not hold, then again we have $\nu_r(N)=\nu_r(q)$, so assume that $r$ divides $n+2$ also. Then a similar argument as above shows that $r$ divides the numerator of $r_3$ if and only if $r$ divides $2n_2+1=n+4$, but this would imply that $r$ divides $(n+4)-(n+2)=2$, an impossibility since $r>3$. This concludes the proof.\qed 

\medskip
These lemmas allow one to deduce in certain cases that $\rho_{m,n}$ is noncongruence, as we now demonstrate. We make no attempt at completeness here, and choose instead to pick off easily obtained consequences of what was just proven. In particular, the primes 2 and 3 (as always in the theory of modular forms) present special difficulties and for this reason have been avoided. More generally, when the level $N$ is a composite number with a large number of primes occurring to high powers it becomes much more tedious to state and prove results analogous to what appears below. For this reason, we restrict ourselves to minimal models $V(p,q)$ where $p$ and $q$ are prime powers. In what follows, we use the ceiling function notation $\lceil x\rceil$,  which denotes the smallest integer $n$ such that $n\geq x$. 
\begin{thm}\label{thm:pq} 
Assume $\rho_{m,n}$ is the representation \eqr{rhomn} associated to the irreducible $V(p,q)$-module $L_{m,n}$, such that $p=r^a,q=s^b$ are powers of distinct primes $r,s>3$, and let $\alpha := \lceil r^{a-2} \rceil$ and $\beta := \lceil s^{b-2} \rceil$. Then $\rho_{m,n}$ is noncongruence so long as the inequalities
\[
\alpha\leq m\leq p-4, \qquad \beta\leq n\leq q-4
\]
and one of $\alpha<m,\ \beta<n$ hold.
\end{thm}

\proof From \eqr{rj} it follows that the level $N$ of $\rho_{m,n}$ is a divisor of $48pq$, and the current assumptions allow us to conclude from Lemmas \ref{lem:pdiv} and \ref{lem:qdiv} that both $p$ and $q$ divide $N$. Both $p$ and $q$ are odd in this setting, and since $m$ and $n$ are odd by convention \eqr{conventions}, one sees from \eqr{rj} that the numerator of each $r_j$ is divisible by 2 but not by 4. Thus 8 divides $N$ as well. One now checks that Theorem \ref{thm:NW} with $p=r^a$, $q=s^b$ and $\alpha,\beta$ defined as above yields 
$$\dim\rho_{m,n}\geq2\cdot\frac{p-\alpha}{2}\cdot\frac{q-\beta}{2}=\frac{(p-\alpha)(q-\beta)}{2}$$
if $\rho_{m,n}$ is congruence. Since we assume that at least one of $m>\alpha,\ n>\beta$ is true, by Corollary \ref{cor:dim} we have that $\rho_{m,n}$ is noncongruence.\qed

\medskip
There are two more cases to consider.
\begin{cor}\label{cor:q-2}
Suppose that one of the following holds:
\begin{enumerate}
\item[(i)]$p\geq3$ is odd, $m=p-2$, and $n,b,s,q$ satisfy the hypotheses of Theorem \ref{thm:pq} with $n>\beta$.
\item[(ii)]$q\geq3$ is odd, $n=q-2$, and $m,a,r,p$ satisfy the hypotheses of Theorem \ref{thm:pq} with $m>\alpha$.
\end{enumerate} 
Then $\rho_{m,n}$ is noncongruence.
\end{cor}

\proof Assume case $(i)$ holds. Since $m=p-2$, we have from Corollary \ref{cor:dim} that 
\[
\dim\rho_{m,n}=\frac{(p-m)(q-n)}{2}=q-n.
\]
Because $m$ is odd, the proof of Theorem \ref{thm:pq} shows that $8q$ divides the level of $\rho_{m,n}$. Theorem \ref{thm:NW} now implies that if $\rho_{m,n}$ is congruence then 
\[
\dim\rho_{m,n}\geq2\cdot\frac{(q-\beta)}{2}=q-\beta.
\]
Since we are assuming that $n>\beta$, it must be that $\rho_{m,n}$ is noncongruence. The proof of case $(ii)$ follows similarly.\qed 

\medskip
The proof of Corollary \ref{cor:pq-2} may now be given. In the notation of Theorem \ref{thm:pq}, we have in this setting that $\alpha=\beta=1$, and the theorem along with Corollary \ref{cor:q-2} covers all the claimed cases.



\end{document}